\documentclass{amsart}
\title{Rothberger's property in all finite powers
\footnote{\lowercase{{\bf {\uppercase{K}}ey words and phrases:} {\uppercase{E}}llentuck topology, {\uppercase{R}}amsey Theory, {\uppercase{R}}othberger property, forcing, selection principle.\\
{{\bf{\uppercase{S}}ubject {\uppercase{C}}lassification:} {\uppercase{P}}rimary 03{\uppercase{E}}02, 05{\uppercase{C}}55, 05{\uppercase{D}}10,  54{\uppercase{D}}20.}
}}
}
\author{by Marion Scheepers}
\usepackage{amsfonts}
\newtheorem{theorem}{Theorem}
\newtheorem{lemma}[theorem]{Lemma}
\newtheorem{corollary}[theorem]{Corollary}

\newtheorem{definition}{Definition}

\newcommand{\pf}{\bf Proof:}
\newcommand{\epf}{\Box\vspace{0.15in}}
\newcommand{\forces}{\mathrel{\|}\joinrel\mathrel{-}}

\newcommand{\sone}{{\sf S}_1}
\newcommand{\gone}{{\sf G}_1}
\newcommand{\sfin}{{\sf S}_{fin}}

\newcommand{\egp}{{\sf E}}

\newcommand{\naturals}{{\mathbb N}}
\begin{document}
\maketitle
\begin{abstract} 
   A space $X$ has the Rothberger property in all finite powers if, and only if, its collection of $\omega$-covers has Ramseyan properties.
\end{abstract}

For $s\in\,[\mathbb{N}]^{<\aleph_0}$ and for $B\in[\mathbb{N}]^{\aleph_0}$ use $s< B$ to denote that $s=\emptyset$ or $\max(s) < \min(B)$. For $s < B$  define
$
  [s,B] = \{s\cup C\in[\mathbb{N}]^{\aleph_0}: \, s < C\subseteq B\}.
$
The family $\{[s,B]:\, s\subset\mathbb{N} \mbox{ finite and } s < B\in[\mathbb{N}]^{\aleph_0}\}$ forms a basis for a topology on $[\mathbb{N}]^{\aleph_0}$.  This is the \emph{Ellentuck topology} on $[\mathbb{N}]^{\aleph_0}$ and was introduced in \cite{El}. 

Recall that a subset $N$ of a topological space is \emph{nowhere dense} if there is for each nonempty open set $U$ of the space a nonempty open subset $V\subset U$ such that $N\cap V = \emptyset$. And $N$ is said to be \emph{meager} if it is a union of countably many nowhere dense sets. A subset of a topological space is said to have the \emph{Baire property} if it is of the form $(U\setminus M)\bigcup(M\setminus U)$ for some 
open set $U$ and some meager set $M$.
\begin{theorem}[Ellentuck]\label{galvinprikry} For a set $R\subset [\mathbb{N}]^{\aleph_0}$ the following are equivalent:
\begin{enumerate}
\item{$R$ has the Baire property in the Ellentuck topology.}
\item{For each finite set $s\subset\naturals$ and for each infinite set $S\subset \mathbb{N}$ with $s< S$ there is an infinite set $T\subset S$ such that
      either $[s,T]\subset R$, or else $[s,T]\cap R = \emptyset$.}
\end{enumerate}
\end{theorem}

The proof of $(1)\Rightarrow(2)$ is nontrivial but uses only the techniques of Galvin and Prikry \cite{G-P}. Galvin and Prikry proved a precursor of  Theorem \ref{galvinprikry}: If $R$ is a Borel set in the topology inherited from $2^{\naturals}$ via representing sets by their characteristic functions, then $R$ has property (2) in Theorem \ref{galvinprikry}. Silver and Mathias subsequently gave metamathematical proofs that analytic sets (in the $2^{\naturals}$-topology) have this property. Theorem \ref{galvinprikry} at once yields all these prior results. The original papers \cite{El} and \cite{G-P} give a nice overview of these facts, and more.

When a subset $\mathcal{X}$ of $[\naturals]^{\aleph_0}$ inherits the Ellentuck topology from $[\naturals]^{\aleph_0}$, we shall speak of ``$\mathcal{X}$ with the Ellentuck topology". For $A$ an abstract countably infinite set define the Ellentuck topology on $[A]^{\aleph_0}$ by fixing a bijective enumeration $(a_n:n\in\naturals)$ of $A$ and by defining for $s$ and $T$ nonempty subsets of $A$:  
\[
  s< T \mbox{ if: }a_n\in s \mbox{ and }a_m\in T \Rightarrow n < m.
\]
With the relation $s< T$ defined, define the Ellentuck topology on $[A]^{\aleph_0}$ as above. For $B\subseteq A$ and for finite set $s\subseteq A$ we write $B|s$ for $\{a_n\in B: s< \{a_n\}\}$.

For families $\mathcal{A}$ and $\mathcal{B}$ we now define a sequence of statements:
\begin{quote}
 $\egp(\mathcal{A},\mathcal{B})$: For each countably infinite $A\in\mathcal{A}$ and for each set $R\subset [A]^{\aleph_0}\cap\mathcal{B}$ the implication (1)$\Rightarrow$(2) holds, where:
\begin{enumerate}
\item{$R$ has the Baire property in the Ellentuck topology on $[A]^{\aleph_0}\cap\mathcal{B}$.}
\item{For each $S\subset A$ with $S\in\mathcal{A}$ and each finite subset $s$ of $A$, there is an infinite $B\subset S|s$ with $B\in\mathcal{B}$ such that $[s,B]\cap\mathcal{B}\subseteq R$ or $[s,B]\cap\mathcal{B}\cap R = \emptyset$.}
\end{enumerate}
\end{quote}
Thus, ${\sf E}([\naturals]^{\aleph_0},[\naturals]^{\aleph_0})$ is Ellentuck's theorem.

\begin{quote}
 ${\sf GP}(\mathcal{A},\mathcal{B})$: For each countably infinite $A\in\mathcal{A}$ and each $R\subset [A]^{\aleph_0}\cap \mathcal{B}$ the implication $(1)\Rightarrow(2)$ holds:
\begin{enumerate}
\item{$R$ is open in the $2^{\naturals}$ topology on $[A]^{\aleph_0}\cap\mathcal{B}$.}
\item{For each $S\in[A]^{\aleph_0}\cap\mathcal{A}$ there is a set $B\in[S]^{\aleph_0}\cap\mathcal{B}$ such that either $([B]^{\aleph_0}\cap\mathcal{B})\subseteq R$, or else $[B]^{\aleph_0}\cap\mathcal{B}\cap R = \emptyset$.} 
\end{enumerate}
\end{quote}
Thus, ${\sf GP}([\naturals]^{\aleph_0},[\naturals]^{\aleph_0})$ is part of the Galvin-Prikry theorem.

\begin{definition} A subset $\mathcal{S}$ of $[A]^{<\aleph_0}$ is:
\begin{enumerate}
\item{dense if for each $B\in[A]^{\aleph_0}\cap \mathcal{A}$, $\mathcal{S}\cap[B]^{<\aleph_0} \neq \emptyset$.}
\item{thin if no element of $\mathcal{S}$ is an initial segment of another element of $\mathcal{S}$.}
\end{enumerate}
\end{definition}
The following is an abstract formulation of Galvin's generalization of Ramsey's Theorem, announced in \cite{Galvinnotices} and in \cite{G-P} derived from Theorem 1 there:
\begin{quote} ${\sf FG}(\mathcal{A},\mathcal{B})$: 
For each countably infinite $A\in\mathcal{A}$ and for each dense set $\mathcal{S}\subset [A]^{<\aleph_0}$ there is a $B\in[A]^{\aleph_0}\cap\mathcal{B}$ such that each $C\in[B]^{\aleph_0}\cap\mathcal{B}$ has an initial segment in $\mathcal{S}$.
\end{quote}
In this notation Galvin's generalization of Ramsey's theorem reads that ${\sf FG}([\naturals]^{\aleph_0},[\naturals]^{\aleph_0})$.
Similarly, the following is an abstract formulation of Nash-Williams' theorem:
\begin{quote}${\sf NW}(\mathcal{A},\mathcal{B})$:
For each countably infinite $A\in\mathcal{A}$ and for each thin family $\mathcal{T}\subset [A]^{<\aleph_0}$ and for each $n$, and each partition $\mathcal{T} = \mathcal{T}_1 \cup \mathcal{T}_2 \cup \cdots \cup \mathcal{T}_n$ there is a $B\in[A]^{\aleph_0}\cap \mathcal{B}$ and an $i\in\{1,\cdots,n\}$ such that $[B]^{<\aleph_0}\cap\mathcal{T} \subseteq \mathcal{T}_i$.
\end{quote} 
In this notation Nash-Williams' theorem reads that ${\sf NW}([\naturals]^{\aleph_0},[\naturals]^{\aleph_0})$.
\begin{quote} 
$\mathcal{A}\longrightarrow(\mathcal{B})^n_k$: For positive integers $n$ and $k$ and for each countable $A\in\mathcal{A}$ and for each function 
$f:[A]^n\rightarrow\{1,\cdots,\, k\}$ there is a $B\in[A]^{\aleph_0}\cap\mathcal{B}$ and an $i\in\{1,\cdots,k\}$ such that $f$ has value $i$ on $[B]^n$.
\end{quote}
In this notation Ramsey's theorem reads: For each $n$ and $k$,  $[\naturals]^{\aleph_0}\longrightarrow([\naturals]^{\aleph_0})^n_k$.

An open cover $\mathcal{U}$ of a topological space $X$ is said to be an $\omega$-cover if $X\not\in\mathcal{U}$, but there is for each finite set $F\subset X$ a $U\in\mathcal{U}$ with $F\subseteq U$. The symbol $\Omega_X$ denotes the collection of $\omega$-covers of $X$. The symbol $\mathcal{O}_X$ denotes the collection of open covers of $X$. In \cite{Ro38} Rothberger introduced the following covering property: For each sequence $(\mathcal{U}_n:n\in\naturals)$ of open covers of $X$ there is a sequence $(U_n:n\in\naturals)$ such that each $U_n\in\mathcal{U}_n$, and $\{U_n:n\in\naturals\}$ is a cover of $X$. The symbol $\sone(\mathcal{O}_X,\mathcal{O}_X)$ denotes this statement. The corresponding statement for $\omega$-covers of $X$, $\sone(\Omega_X,\Omega_X)$, was introduced in \cite{Sakai} by Sakai. It states: For each sequence $(\mathcal{U}_n:n\in\naturals)$ of $\omega$-covers of $X$ there is a sequence $(U_n:n\in\naturals)$ such that each $U_n\in\mathcal{U}_n$, and $\{U_n:n\in\naturals\}$ is an $\omega$ cover for $X$. Sakai proved that $X$ has $\sone(\Omega_X,\Omega_X)$ if, and only if, all finite powers of $X$ have $\sone(\mathcal{O}_X,\mathcal{O}_X)$. According to Gerlits and Nagy \cite{GN} a space is said to be an $\epsilon$-\emph{space} if each $\omega$-cover contains a countable subset which still is an $\omega$-cover. A space is an $\epsilon$-space if and only if it has the Lindel\"of property in all finite powers - see \cite{GN} for details. In this paper we prove:

\begin{theorem}\label{egprothberger} For an $\epsilon$-space $X$, the following are equivalent:
\begin{enumerate}
\item{$\sone(\Omega_X,\Omega_X)$.}
\item{$\egp(\Omega_X,\Omega_X)$.}
\item{${\sf GP}(\Omega_X,\Omega_X)$.}
\item{${\sf FG}(\Omega_X,\Omega_X)$.}
\item{${\sf NW}(\Omega_X,\Omega_X)$.}
\item{For all $n$ and $k$, $\Omega_X\rightarrow(\Omega_X)^n_k$.}
\item{$\Omega_X\rightarrow(\Omega_X)^2_2$.}
\end{enumerate}
\end{theorem}

\section{The proof of $\sone(\Omega_X,\Omega_X)\,\Rightarrow\,\egp(\Omega_X,\Omega_X)$:}

Assume that $X$ has property $\sone(\Omega_X,\Omega_X)$. Fix a countable $A\in\Omega_X$ and fix a set $R\subset [A]^{\aleph_0}\cap\Omega_X$. %By $\sone(\Omega_X,\Omega_X)$ we may assume that $A$ is countable. 
For the remainder of the argument, fix a bijective enumeration of $A$, say $(a_n:n\in\naturals)$. 
Sets of the form $[s,C] = \{D: s<C \mbox{ and }s \subset D \subseteq s\cup C\}$ constitute a basis for the Ellentuck topology on $[A]^{\aleph_0}$.

\begin{definition} For a finite set $s\subset A$ and for $B\in[A|s]^{\aleph_0}\cap\Omega_X$:
\begin{enumerate}
\item{$B$ accepts $s$ if $[s,B]\cap\Omega_X\subseteq R$.}
\item{$B$ rejects $s$ if no $C\in[B]^{\aleph_0}\cap\Omega_X$ accepts $s$.}
\end{enumerate}
\end{definition}
Lemma \ref{acceptreject} will be used without special reference:
\begin{lemma}\label{acceptreject} Let a finite set $s\subset A$ and a set $B\in[A|s]^{\aleph_0}\cap\Omega_X$ be given:
\begin{enumerate}
\item{$B$ accepts $s$ if, and only if, each $C\in[B]^{\aleph_0}\cap\Omega_X$ accepts $s$.}
\item{$B$ rejects $s$ if, and only if, each $C\in[B]^{\aleph_0}\cap\Omega_X$ rejects $s$.}
\end{enumerate}
\end{lemma}

\begin{lemma}\label{decideoneset} For each finite set $s\subset A$, there is a $B\in[A|s]^{\aleph_0}\cap \Omega_X$ such that $B$ accepts $s$ or $B$ rejects $s$.
\end{lemma}
$\pf$ If $A|s$ does not reject $s$, choose a $B\in[A|s]^{\aleph_0}\cap\Omega_X$ accepting $s$. $\epf$

\begin{lemma}\label{acceptstrongreject} Let $t\subset A$ be a finite set. Let $B\in[A]^{\aleph_0}\cap\Omega_X$ be such that for each finite set $s\subset (t\cup B)$, $B|s$ accepts $s$ or $B|s$ rejects $s$. If $B|t$ rejects $t$ then $C =\{u\in B: B|(t\cup\{u\}) \mbox{ rejects } t\cup\{u\}\}$ is a member of $\Omega_X$. 
\end{lemma}
$\pf$
Suppose not. Then $D = t\cup(B\setminus C)\in\Omega_X$, and for each $u\in D|t$, $B|(t\cup\{u\})$ accepts $t\cup\{u\}$. Thus for each $u\in D|t$, $D|(t\cup\{u\})$ accepts $t\cup\{u\}$. This means that $[t,D|t] = \cup_{u\in D}[t\cup\{u\},D|(t\cup\{u\})]\subseteq R$, and so $D|t$ accepts $t$. This contradicts Lemma \ref{acceptreject} (2) since $D\in [B]^{\aleph_0}\cap\Omega_X$ and $B|t$ rejects $t$.
$\epf$

\section*{$\omega$-covers accepting or rejecting all finite subsets.}

The game $\gone(\Omega_X,\Omega_X)$ is played as follows: Players ONE and TWO play an inning per positive integer. In the $n$-th inning ONE first chooses an $O_n\in\Omega_X$; TWO responds with a $T_n\in O_n$. A play $O_1,\, T_1,\, \cdots,\, O_n,\, T_n,\, \cdots$ is won by TWO if $\{T_n:n\in\naturals\}\in \Omega_X$; else, ONE wins. It was shown in \cite{coc3} that

\begin{theorem}\label{coc3gone} For a topological space $Y$ the following are equivalent:
\begin{enumerate}
\item{$Y$ has property $\sone(\Omega_Y,\Omega_Y)$.}
\item{ONE has no winning strategy in $\gone(\Omega_Y,\Omega_Y)$.}
\end{enumerate}
\end{theorem}

\begin{theorem}\label{decidedfin} If $Y$ has property $\sone(\Omega_Y,\Omega_Y)$, then for each finite set $t\subset A$ and for each $B\in[A|t]^{\aleph_0}\cap\Omega_Y$ there is a $C\in[B]^{\aleph_0}\cap\Omega_Y$ such that for each finite set $s\subset t\cup C$, $C|s$ accepts $s$ or $C|s$ rejects $s$.
\end{theorem}
$\pf$ Let $t$ and $B\in[A|t]^{\aleph_0}\cap\Omega_Y$ be given. Define a strategy $\sigma$ for ONE of $\gone(\Omega_Y,\Omega_Y)$ as follows:

Enumerate the set of all subsets of $t$ as $\{t_1,\cdots,t_n\}$. Using Lemma \ref{decideoneset} recursively choose $B_1\supset B_2\supset\cdots \supset B_n$ in $[B]^{\aleph_0}\cap \Omega_Y$ such that for each $i$, $B_i$ accepts $t_i$ or $B_i$ rejects $t_i$. Then define:
\[
  \sigma(\emptyset) = B_n.
\]

If TWO now chooses $T_1\in\sigma(\emptyset)$ then use Lemma \ref{decideoneset} in the same way to choose 
\[
  \sigma(T_1) \in[\sigma(\emptyset)|\{T_1\}]^{\aleph_0}\cap\Omega_Y
\]
such that for each set $F\subset t\cup\{T_1\}$, $\sigma(T_1)$ accepts $F$, or rejects $F$.

When TWO responds with $T_2\in F(T_1)$, enumerate the subsets of $t\cup \{T_1, T_2\}$ as $(t_1,\cdots,t_n)$ say, and choose by Lemma \ref{decideoneset} sets $B_1,\cdots B_n\in[\sigma(T_1)|\{T_2\}]^{\aleph_0}\cap \Omega_X$ such that $B_j$ accepts $t_j$ or $B_j$ rejects $t_j$ for $1\le j\le n$ and $B_j\subset B_{j-1}$. Finally put 
\[
  \sigma(T_1,T_2) = B_n.
\]
Note that for each finite subset $F$ of $t\cup\{T_1, T_2\}$, $\sigma(T_1,T_2)$ accepts $F$ or rejects it.

It is clear how player ONE's strategy is defined. By Theorem \ref{coc3gone} $\sigma$ is not a winning strategy for ONE. Consider a $\sigma$-play lost by ONE, say
\[
  \sigma(\emptyset), \, T_1,\, \sigma(T_1),\, T_2,\, \sigma(T_1,T_2),\, \cdots,\, T_n,\, \sigma(T_1,\cdots,T_n),\, \cdots
\]
Then $C = t\cup \{T_n:\, n\in\naturals\} \subset B$ is an element of $\Omega_Y$.

We claim that for each finite subset $s$ of $t\cup C$, $C|s$ accepts $s$ or $C|s$ rejects $s$. For consider such a $s$. If $s\subseteq t$, then as $C\subset F(\emptyset)$ and $F(\emptyset)$ accepts or rejects $s$, also $C$ does. If $s\not\subseteq t$, then put $n = \max\{m:T_m\in s\}$. Then $s$ is a subset of $t\cup\{T_1,\cdots,T_n\}$, so that $s$ is accepted or rejected by $\sigma(T_1,\cdots,T_n)$. But $C|s \subseteq \sigma(T_1,\cdots,T_n)$, and so $C|s$ accepts or rejects $s$. 
$\epf$

\section*{Completely Ramsey sets}
 
The subset $R$ of $[A]^{\aleph_0}\cap \Omega_X$ is said to be \emph{completely Ramsey} if there is for each finite set $s\subset A$ and for each $B\in[A|s]^{\aleph_0}\cap\Omega_X$ a set $C\in[B]^{\aleph_0}\cap\Omega_X$ such that 
\begin{enumerate}
\item{either $([s,C]\cap\Omega_X)\subseteq R$,}
\item{or else $([s,C]\cap\Omega_X)\cap R = \emptyset$.}
\end{enumerate}

\begin{lemma}\label{complramseyunions} If $R$ and $S$ are completely Ramsey subsets of $[A]^{\aleph_0}\cap \Omega_X$, then so is $R\bigcup S$.
\end{lemma}
$\pf$ Let a finite set $s\subset A$ and $B\in[A|s]^{\aleph_0}\cap\Omega_X$ be given. Since $R$ is completely Ramsey, choose $C\in[B]^{\aleph_0}\cap\Omega_X$ such that $([s,C]\cap\Omega_X)\subset R$, or $([s,C]\cap\Omega_X)\cap R = \emptyset$. If the former hold we are done. In the latter case, since $S$ is completely Ramsey, choose $D\in[C]^{\aleph_0}\cap\Omega_X$ such that $([s,D]\cap\Omega_X)\subseteq S$, or $([s,D]\cap\Omega_X)\cap S = \emptyset$. In either case the proof is complete. $\epf$

The following Lemma is obviously true.
\begin{lemma}\label{complramseycomplements} If $R$ is completely Ramsey, then so is $([A]^{\aleph_0}\cap\Omega_X)\setminus R$.
\end{lemma}

\begin{corollary}\label{intersections} If $R$ and $S$ are completely Ramsey subsets of $[A]^{\aleph_0}\cap \Omega_X$, then so is $R\bigcap S$.
\end{corollary}
$\pf$ Lemmas \ref{complramseyunions} and \ref{complramseycomplements}, and De Morgan's laws.$\epf$

\section*{Open sets in the Ellentuck topology}

We are still subject to the hypothesis that $X$ satisfies $\sone(\Omega_X,\Omega_X)$. 
\begin{lemma}\label{openprecursor} For each finite set $t\subset A$ and for each  $B\in[A|t]^{\aleph_0}\cap\Omega_X$ such that for each finite subset $F$ of $t\cup B$,  $B|F$ accepts, or rejects $F$ the following holds: For each finite set $s\subset t\cup B$ such that $B|s$ rejects $s$, there is a $C\in[B|s]^{\aleph_0}\cap\Omega_X$ such that for each finite set $F\subset C$, $C|F$ rejects $s\cup F$.
\end{lemma}
$\pf$ Fix $B$ and $s$ as in the hypotheses. Define a strategy $\sigma$ for ONE in $\gone(\Omega_X,\Omega_X)$ as follows: 
By Lemma \ref{acceptstrongreject} 
\[
  \sigma(\emptyset) =\{U\in B:\,s<\{U\}\mbox{ and $B|\{U\}$ rejects }s\cup\{U\}\} \in\Omega_X.
\]
 Notice that $\sigma(\emptyset)$ accepts or rejects each of its finite subsets, it rejects $s$, and for each $U\in \sigma(\emptyset)$, $\sigma(\emptyset)|\{U\}$ rejects $s\cup\{U\}$.

If TWO now chooses $T_1\in \sigma(\emptyset)$, then by  Lemma \ref{acceptstrongreject} 
\[
  \sigma(T_1)=\{U\in \sigma(\emptyset)\setminus\{T_1\}: \sigma(\emptyset)|F \mbox{ rejects }s\cup F \mbox{ for each finite }F\subset\{T_1,U\}\}
\]
 is in $\Omega_X$. As before, $\sigma(T_1)$ accepts or rejects each of its finite subsets, and for any $U\in \sigma(T_1)$, for each finite subset $F$ of $\{U,T_1\}$, $\sigma(T_1)|F$ rejects $s\cup F$.  

If next TWO chooses $T_2\in \sigma(T_1)$, then by Lemma \ref{acceptstrongreject}
\[
  \sigma(T_1,T_2) = \{U\in \sigma(T_1)\setminus\{T_2\}: \sigma(T_1)|F \mbox{ rejects } s\cup F \mbox{ for any finite } F\subset \{T_1,T_2,U\}\}
\]
is an element of $\Omega_X$.

Continuing in this way we define a strategy $\sigma$ for ONE in $\gone(\Omega_X,\Omega_X)$. Since $X$ satisfies $\sone(\Omega_X,\Omega_X)$, $\sigma$ is not a winning strategy for ONE. Consider a $\sigma$-play lost by ONE, say:
\[
 \sigma(\emptyset), \, T_1,\, \sigma(T_1),\, T_2,\, \sigma(T_1,T_2),\, T_3,\, \sigma(T_1,T_2,T_3),\, \cdots
\]
Put $C = \{T_n:n\in\naturals\}$. Then $C\in [B|s]^{\aleph_0}\cap\Omega_X$. We claim that for each finite set $F\subset C$, $C|F$ rejects $s\cup F$. 

For choose a finite set $F\subset C$. Then $F\cap s = \emptyset$. Fix $n=\max\{m:T_m\in F\}$. Then $C|F \subset \sigma(T_1,\cdots,T_n)$, and the latter rejects $s\cup F$ for all finite subsets $F$ of $\{T_1,\cdots,T_n\}$. Thus $C|F$ rejects $s\cup F$.$\epf$ 

\begin{theorem}\label{opencomplRamsey} If $X$ has property $\sone(\Omega_X,\Omega_X)$, then every open subset of $[A]^{\aleph_0}\cap\Omega_X$ is completely Ramsey.
\end{theorem}
$\pf$ Let $R\subset[A]^{\aleph_0}\cap\Omega_X$ be open in this subspace. Consider a finite set $s\subset A$ and a $B \in [A|s]^{\aleph_0}\cap\Omega_X$. Since $(X,d)\models\sone(\Omega_X,\Omega_X)$, choose by Theorem \ref{decidedfin} a $C\in[B]^{\aleph_0}\cap\Omega_X$ such that for each finite set $F\subset(s\cup C)$, $C|F$ accepts or rejects $F$.

If $C$ accepts $s$ then we have $[s,C]\cap\Omega_X\subseteq R$, and we are done. Thus, assume that $C$ does not accept $s$. Then $C$ rejects $s$, and we choose by Lemma \ref{openprecursor} a $D\in[C|s]^{\aleph_0}\cap\Omega_X$ such that for each finite subset $F$ of $D$, $D|F$ rejects $s\cup F$.

We claim that $([s,D]\cap\Omega_X)\cap R = \emptyset$. For suppose not. Choose $E\in[s,D] \cap\Omega_X \cap R$. Since $R$ is open, choose an Ellentuck neighborhood of $E$ contained in $R$, say $[t,K]\cap \Omega_X$. Then we have $s\subset E \subset s\cup D$ and $t\subset E \subset t\cup K$. But then $s\cup t\subset E \subset t \cup K$ and $[s\cup t,K|s]\subset R$, whence also $[s\cup t,E|(s\cup t)]\subset R$. But then $E|(s\cup t)$ accepts $s\cup t$ where $t$ is a finite subset of $s\cup D$, and $E|(s\cup t)\subset D|t$, and $D|t$ rejects $s\cup t$, a contradiction. $\epf$

\section*{Meager subsets in the Ellentuck topology}

If the subset $R$ of $[A]^{\aleph_0}\cap\Omega_X$ is nowhere dense in the topology, then for each $B\in[A]^{\aleph_0}\cap\Omega_X$ and for each finite set $s\subset A$, $B|s$ rejects $s$. We now examine the meager subsets of $[A]^{\aleph_0}\cap \Omega_X$.

\begin{lemma}\label{nwdfin} If $R$ is nowhere dense, then there is for each $B\in[A]^{\aleph_0}\cap\Omega_X$ and each finite set $t\subset A$ a set $C\in[B|t]^{\aleph_0}\cap\Omega_X$ such that for each finite set $s\subset t\cup C$, $C|s$ rejects $s$.
\end{lemma}
$\pf$ Since $R$ is nowhere dense, no $\omega$-cover contained in $A$ can accept a finite set. Thus each $\omega$-cover contained in $A$ rejects each finite subset of $A$. $\epf$

\begin{lemma}\label{clnwdcomplramsey} Assume $\sone(\Omega_X,\Omega_X)$. 
If $R$ is a closed nowhere dense subset of $[A]^{\aleph_0}\cap\Omega_X$ then there is for each finite subset $s\subset A$ and for each $B\in[A|s]^{\aleph_0}\cap\Omega_X$ a $C\in[B]^{\aleph_0}\cap\Omega_X$ such that $[s,C]\cap R = \emptyset$. 
\end{lemma}
$\pf$ First, note that closed nowhere dense subsets are complements of open dense sets. By Theorem \ref{opencomplRamsey}, each open set is completely Ramsey. By Lemma \ref{complramseycomplements} each closed, nowhere dense set is completely Ramsey. By Lemma \ref{nwdfin} the rest of the statement follows.
$\epf$

By taking closures, the preceding lemma implies:
\begin{corollary}\label{nwdcomplramsey} Assume $\sone(\Omega_X,\Omega_X)$. If $R$ is a nowhere dense subset of $[A]^{\aleph_0}\cap\Omega_X$ then there is for each finite subset $s\subset A$ and for each $B\in[A|s]^{\aleph_0}\cap\Omega_X$ a $C\in[B]^{\aleph_0}\cap\Omega_X$ such that $[s,C]\cap R = \emptyset$. 
\end{corollary}

And now we prove:
\begin{theorem}\label{meagerellentuck} Assume $\sone(\Omega_X,\Omega_X)$. For a subset $N$ of $[A]^{\aleph_0}\cap\Omega_X$ the following are equivalent:
\begin{enumerate}
\item{$N$ is nowhere dense.}
\item{$N$ is meager.}
\end{enumerate}
\end{theorem}
$\pf$ We must show that (2)$\Rightarrow$(1). Thus, assume that $N$ is meager and write $N=\bigcup_{n\in\naturals}N_n$, where for each $n$ we have $N_n\subseteq N_{n+1}$, and $N_n$ is nowhere dense in $[A]^{\aleph_0}\cap \Omega_X$. Consider any basic open set $[s,B]$ of $[A]^{\aleph_0}\cap\Omega_X$. 
Define a strategy $\sigma$ for ONE in the game $\gone(\Omega_X,\Omega_X)$ as follows: 

Since $N_1$ is nowhere dense, choose by Corollary \ref{nwdcomplramsey} an $O_1\in[B]^{\aleph_0}\cap \Omega_X$ with $[s,O_1]\cap N_1 = \emptyset$. Define $\sigma(\emptyset) = O_1$. 

When TWO chooses $T_1\in \sigma(\emptyset)$ choose by Corollary \ref{nwdcomplramsey} an $O_2\in[\sigma(\emptyset)|\{T_1\}]^{\aleph_0}\cap\Omega_X$ with $[s,O_2]\cap N_2 = \emptyset$, and define $\sigma(T_1) = O_2$.

Now when TWO chooses $T_2\in \sigma(T_1)$, find by Corollary \ref{nwdcomplramsey} an $O_3\in[\sigma(T_1)|\{T_2\}]^{\aleph_0}\cap\Omega_X$ with $[s,O_3]\cap N_3 = \emptyset$, and define $\sigma(T_1,T_2) = O_3$.

It is clear how to define ONE's strategy $\sigma$. By Theorem \ref{coc3gone} $F$ is not a winning strategy for ONE. Consider a play
\[
  \sigma(\emptyset),\, T_1,\, \sigma(T_1),\, T_2,\, \cdots,\, \sigma(T_1,\cdots,T_n),\, T_{n+1},\, \cdots
\] 
lost by ONE. Put $C = \{T_n:n\in\naturals\}$. Then $C\in[B]^{\aleph_0}\cap\Omega_X$. Observe that by the definition of $\sigma$ we have for each $k$ and each finite set $F\subset \{T_1,\cdots,T_k\}$ that $[s\cup F,\sigma(T_1,\cdots,T_k)]\cap N_k = \emptyset$.\\
{\flushleft{\bf Claim:}} $[s,C]\cap N = \emptyset$. \\
For suppose that instead $[s,C]\cap N \neq \emptyset$. Choose $V\in [s,C]\cap N$, and then choose $m$ so that $V\in N_m$. Choose the least $k>m$ with $T_k\in V| s$. This is possible because $s$ is finite. Observe also that $s\subseteq V\subseteq s\cup C = s\cup\{T_j:j\in\naturals\}$. Put $F = V\cap \{T_1,\cdots,T_k\}$. Thus we have that $[s\cup F,V|F]\cap N_k \neq \emptyset$, which contradicts the fact that  $V|F\subset \sigma(T_1,\cdots,T_k)$, and $[s\cup F,\sigma(T_1,\cdots,T_k)]\cap N_k = \emptyset$. This completes the proof of the claim.
$\epf$

Using Lemmas \ref{complramseyunions} and \ref{complramseycomplements} and Corollary \ref{intersections} we have:
\begin{theorem}\label{ellentuck} Suppose $X$ satisfies $\sone(\Omega_X,\Omega_X)$. Then for each $A\in\Omega_X$, every subset of $[A]^{\aleph_0}\cap\Omega_X$ which has the Baire property is completely Ramsey.
\end{theorem}

\section{The proof of $\egp(\Omega_X,\Omega_X)\,\Rightarrow\,\sone(\Omega_X,\Omega_X)$:}

Note that a set open in the $2^{\naturals}$ topology is also open in the Ellentuck topology. The implication $(2)\Rightarrow(3)$ of Theorem \ref{egprothberger} follows from this remark. Now we start with $(3)$. %Define: A subset $R$ of $[A]^{\aleph_0}\cap \Omega_X$ is said to be \emph{Ramsey} if there is for each $B\in[A]^{\aleph_0}\cap\Omega_X$ a set $C\in[B]^{\aleph_0}\cap\Omega_X$ such that either $([C]^{\aleph_0}\cap\Omega_X)\subseteq R$, or else $[C]^{\aleph_0}\cap\Omega_X\cap R = \emptyset$. 

\begin{lemma}\label{densetoinitsegm} Assume ${\sf GP}(\Omega_X,\Omega_X)$. Then ${\sf FG}(\Omega_X,\Omega_X)$ holds.
\end{lemma}
$\pf$ Let $\mathcal{S}\subset[A]^{<\aleph_0}$ be dense and define $\mathcal{I}$ to be the set $\{D\in[A]^{\aleph_0}\cap\Omega_X:\, D \mbox{ has an initial segment in }\mathcal{S}\}$. Then we have: 
\[
  \mathcal{I} = \cup\{[s,D|s]:s\in \mathcal{S},\, D\in[A]^{\aleph_0}\cap\Omega_X \mbox{ and s an initial segment of }D\}
\] 
is a $2^{\naturals}$-open subset of $[A]^{\aleph_0}\cap\Omega_X$. Choose a $B\in[A]^{\aleph_0}\cap\Omega_X$ such that $[B]^{\aleph_0}\cap\Omega_X\subset\mathcal{I}$, or $[B]^{\aleph_0}\cap\Omega_X\cap\mathcal{I} = \emptyset$. But the second alternative implies the contradiction that $[B]^{<\aleph_0}\cap\mathcal{S} = \emptyset$.  It follows that the first alternative holds.$\epf$

\begin{theorem}\label{thintohomogeneous} Assume ${\sf FG}(\Omega_X,\Omega_X)$. Then ${\sf NW}(\Omega_X,\Omega_X)$ holds. 
\end{theorem}
$\pf$ Fix a thin family $\mathcal{T}\subset [A]^{<\aleph_0}$ and positive integer $n$, and a partition $\mathcal{T} = \mathcal{T}_1 \cup \mathcal{T}_2 \cup \cdots \cup \mathcal{T}_n$. 
We may assume $n=2$. If $\mathcal{T}_1$ is not dense, we can choose $B\in[A]^{\aleph_0}\cap \Omega_X$ such that $[B]^{<\aleph_0}\cap\mathcal{T} \subseteq \mathcal{T}_2$. Thus, assume $\mathcal{T}_1$ is dense. Choose, by the hypothesis, 
a $B\in[A]^{\aleph_0}\cap\Omega_X$ such that for each $C\in[B]^{\aleph_0}\cap\Omega_X$, some initial segment of $C$ is in $\mathcal{T}_1$.

Consider any $s\in\mathcal{T}\cap[B]^{<\aleph_0}$, and put $D = s\cup (B|s)$. Then $s$ is an initial segment of $D$, and $D\in[B]^{\aleph_0}\cap\Omega_X$, and so some initial segment of $D$, say $t$, is in $\mathcal{T}_1$. Since both $t$ and $s$ are initial segments of $D$ and are both in $\mathcal{T}$, and since $\mathcal{T}$ is thin, we have $s=t$, and so $s\in\mathcal{T}_1$. Consequently we have $[B]^{<\aleph_0}\cap\mathcal{T}\subseteq \mathcal{T}_1$.
$\epf$

\begin{theorem}\label{ramsey}  Assume that ${\sf NW}(\Omega_X,\Omega_X)$ holds.  
Then: For each $n$ and $k$ we have $\Omega_X\rightarrow(\Omega_X)^n_k$.
\end{theorem}
$\pf$ Let $A\in\Omega_X$ be countable. Let positive integers $n$ and $k$ be given. Put $\mathcal{T}=[A]^n$. Then $\mathcal{T}$ is thin. Apply the hypothesis.
$\epf$

The following theorem was proven in \cite{coc2} (Theorem 6.1) and \cite{coc1} (Theorem 24)\footnote{See Appendix A}. It, together with the above sequence of implications, completes the proof of Theorem \ref{egprothberger}.
\begin{theorem}\label{soneramsey} The following are equivalent:
\begin{enumerate}
\item{For each $n$ and $k$, $\Omega_X\longrightarrow(\Omega_X)^n_k$}
\item{$X\models\sone(\Omega_X,\Omega_X)$.}
\end{enumerate}
\end{theorem}

\section{Remarks}

The results above are given for $\Omega$, but a study of the proofs will reveal that these equivalences hold for several other families $\mathcal{A}$. The main requirements on $\mathcal{A}$ are that each element of $\mathcal{A}$ has a countable subset in $\mathcal{A}$, that for each $k$ $\mathcal{A}\rightarrow (\mathcal{A})^1_k$ holds, and that $\sone(\mathcal{A},\mathcal{A})$ is equivalent to ONE not having a winning strategy in $\gone(\mathcal{A},\mathcal{A})$, and that this is equivalent to $\mathcal{A}\rightarrow(\mathcal{A})^2_2$. Though this general treatment can be given without much additional effort, I preferred to illustrate the equivalences using a well-known concrete example, because of the connections of this example with forcing (pointed out below) and with the famous Borel Conjecture. Here are a few examples of such families $\mathcal{A}$:

For a topological space $X$ and an element $x\in X$, define $\Omega_x = \{A\subset X\setminus\{x\}:\, x\in\overline{A}\}$. According to \cite{Sakai} $X$ has  
strong countable fan tightness at $x$ if the selection principle $\sone(\Omega_x,\Omega_x)$ holds. Consider for a Tychonoff space $X$ the subspace of the Tychonoff product $\Pi_{x\in X}\mathbb{R}$ consisting of the continuous functions from $X$ to $\mathbb{R}$. The symbol ${\sf C}_p(X)$ denotes this subspace with the inherited topology. Since ${\sf C}_p(X)$ is homogeneous, the truth of $\sone(\Omega_f,\Omega_f)$ at some point $f$ implies the truth of $\sone(\Omega_f,\Omega_f)$ at any point $f$. Thus we may confine attention to $\Omega_{\bf o}$, where ${\bf o}$ is the function which is zero on $X$. Using the techniques above one can prove:
\begin{theorem}\label{strfantight} For a Tychonoff space $X$ the following are equivalent for ${\sf C}_p(X)$:
\begin{enumerate}
\item{$\sone(\Omega_{\bf o},\Omega_{\bf o})$.}
\item{$\egp(\Omega_{\bf o},\Omega_{\bf o})$.}
\item{${\sf GP}(\Omega_{\bf o},\Omega_{\bf o})$.}
\item{${\sf FG}(\Omega_{\bf o},\Omega_{\bf o})$.}
\item{${\sf NW}(\Omega_{\bf o},\Omega_{\bf o})$.}
\item{For all $n$ and $k$, $\Omega_{\bf o}\rightarrow(\Omega_{\bf o})^n_k$.}
\end{enumerate}
\end{theorem}

For a topological space $X$ let $\mathcal{D}$ denote the collection whose members are of the form $\mathcal{U}$, a family of open subsets of $X$, such that no element of $\mathcal{U}$ is dense in $X$, but $\cup\mathcal{U}$ is dense in $X$. And let $\mathcal{D}_{\Omega}$ be the set of $\mathcal{U}\in\mathcal{D}$ such that for each finite family $\mathcal{F}$ of nonempty open subsets of $X$ there is a $U\in\mathcal{U}$ with $U\cap F\neq\emptyset$ for each $F\in\mathcal{F}$. The families $\mathcal{D}$ and $\mathcal{D}_{\Omega}$ were considered in \cite{coc5} where it was proved that for $X$ a set of real numbers, and ${\sf PR}(X)$ the Pixley-Roy space over $X$, the following holds:
\begin{theorem}\label{prreals} If $X$ is a set of real numbers, the following are equivalent for ${\sf PR}(X)$:
\begin{enumerate}
\item{$\sone(\mathcal{D}_{\Omega},\mathcal{D}_{\Omega})$.}
\item{ONE has no winning strategy in the game $\gone(\mathcal{D}_{\Omega},\mathcal{D}_{\Omega})$.}
\item{For each $n$ and $k$ $\mathcal{D}_{\Omega}\rightarrow(\mathcal{D}_{\Omega})^n_k$.}
\end{enumerate}
Each of these statements is equivalent to $X$ having $\sone(\Omega_X,\Omega_X)$.
\end{theorem}
Using the techniques above one can prove:
\begin{theorem}\label{prRamsey} For a set $X$ of reals the following are equivalent for ${\sf PR}(X)$:
\begin{enumerate}
\item{$\sone(\mathcal{D}_{\Omega},\mathcal{D}_{\Omega})$.}
\item{$\egp(\mathcal{D}_{\Omega},\mathcal{D}_{\Omega})$.}
\item{${\sf GP}(\mathcal{D}_{\Omega},\mathcal{D}_{\Omega})$.}
\item{${\sf FG}(\mathcal{D}_{\Omega},\mathcal{D}_{\Omega})$.}
\item{${\sf NW}(\mathcal{D}_{\Omega},\mathcal{D}_{\Omega})$.}
\item{For all $n$ and $k$, $\mathcal{D}_{\Omega}\rightarrow(\mathcal{D}_{\Omega})^n_k$.}
\end{enumerate}
\end{theorem}

For a non-compact topological space $X$ call an open cover $\mathcal{U}$ a {\sf k}-cover if there is for each compact $C\subset X$ a $U\in\mathcal{U}$ such that $C\subseteq U$, and if $X\not\in\mathcal{U}$. Let $\mathcal{K}$ denote the collection of {\sf k}-covers of such an $X$. If $X$ is a separable metric space then each member of $\mathcal{K}$ has a countable subset which still is a member of $\mathcal{K}$. Using the techniques above one can prove:
\begin{theorem}\label{kcovRamsey} For separable metric spaces $X$ the following are equivalent:
\begin{enumerate}
\item{ONE has no winning strategy in $\gone(\mathcal{K},\mathcal{K})$.}
\item{$\sone(\mathcal{K},\mathcal{K})$.}
\item{$\egp(\mathcal{K},\mathcal{K})$.}
\item{${\sf GP}(\mathcal{K},\mathcal{K})$.}
\item{${\sf FG}(\mathcal{K},\mathcal{K})$.}
\item{${\sf NW}(\mathcal{K},\mathcal{K})$.}
\item{For all $n$ and $k$, $\mathcal{K}\rightarrow(\mathcal{K})^n_k$.}
\end{enumerate}
\end{theorem}
The equivalence of (2) and (7) for n=2 and k=2 is Theorem 8 of \cite{dkm}. The equivalence of (1) and (2) is a result of \cite{ST}. The remaining equivalences are then derived as was done above for $\Omega$.

A collection $\mathcal{C}$ of subsets of a set $S$ is said to be a \emph{combinatorial} $\omega$-\emph{cover} of $S$ if $S\not\in\mathcal{C}$, but for each finite subset $F$ of $S$ there is a $C\in\mathcal{C}$ with $F\subseteq C$. For an infinite cardinal number $\kappa$ let $\Omega_{\kappa}$ be the set of \emph{countable} combinatorial $\omega$-covers of $\kappa$. Let ${\sf cov}(\mathcal{M})$ be the least infinite cardinal number $\kappa$ such that the real line is a union of $\kappa$ first category sets. By the Baire Category Theorem ${\sf cov}(\mathcal{M})$ is uncountable. 
Using the techniques of this paper one can prove:
\begin{theorem}\label{cardinality} For an infinite cardinal number $\kappa$ the following are equivalent:
\begin{enumerate}
\item{$\kappa<{\sf cov}(\mathcal{M})$.}
\item{$\sone(\Omega_{\kappa},\Omega_{\kappa})$.}
\item{$\egp(\Omega_{\kappa},\Omega_{\kappa})$.}
\item{${\sf GP}(\Omega_{\kappa},\Omega_{\kappa})$.}
\item{${\sf FG}(\Omega_{\kappa},\Omega_{\kappa})$.}
\item{${\sf NW}(\Omega_{\kappa},\Omega_{\kappa})$.}
\item{For all positive integers $n$ and $k$, $\Omega_{\kappa}\rightarrow(\Omega_{\kappa})^n_k$.}
\end{enumerate}
\end{theorem}

\section{Rothberger's property and forcing}

Now we explore the connections between forcing and Rothberger's property. Much of this part of the paper is inspired by Theorem 9.3 of \cite{JEB}.  

We begin by defining the following version of the well-known \emph{Mathias reals} partially ordered set. Fix as before a countable $\omega$-cover $A$ of $X$, and enumerate it bijectively as $(a_n:n\in\naturals)$. For $s\subset A$ finite, and $C\subset A|s$ with $C\in\Omega_X$, define: 
\[
  \mathcal{M}_A :=\{(s,C): s\in[A]^{<\aleph_0} \mbox{ and }C\subset A|s \mbox{ and }C\in\Omega_X\}.
\]
For $(s_1,C_1)$ and $(s_2,C_2)$ elements of $\mathcal{M}_A$, we define $(s_1,C_1)\prec (s_2,C_2)$ if: $s_2\subset s_1$ and $C_1\subset C_2$ and $s_1\setminus s_2\subset C_2|s_2$. 

Now $(\mathcal{M}_A,\prec)$ is a partially ordered set. Its combinatorial and forcing properties are related to the combinatorial properties of $\omega$-covers of $X$. In this section we will show (see Theorem 9.3 of \cite{JEB}):
\begin{theorem}\label{truthlemma} The following are equivalent:
\begin{enumerate}
\item{$\sone(\Omega_X,\Omega_X)$ holds.}
\item{For each countable $A\in\Omega_X$, for each sentence $\psi$ in the $\mathcal{M}_A$-forcing language, and for each $(s,B)\in\mathcal{M}_A$, there is a $C\subset B$ with $C\in\Omega_X$ such that $(s,C)\forces \psi$, or $(s,C)\forces \neg\psi$}
\end{enumerate}
\end{theorem}

\subsection{Proof of $(1)\Rightarrow (2)$:} Fix a sentence $\psi$ of the $\mathcal{M}_A$-forcing language and fix $(s,B)\in \mathcal{M}_A$. Define the subsets
\[
  \mathcal{W} = \cup\{[t,C]: (t,C)\in\mathcal{M}_A \mbox{ and } (t,C)\forces \psi\}
\]
and
\[
  \mathcal{D} = \cup\{[t,C]: (t,C)\in\mathcal{M}_A \mbox{ and } (t,C)\forces \neg\psi\}.
\]
Then $\mathcal{W}$ and $\mathcal{D}$ are open sets in the Ellentuck topology on $[A]^{\aleph_0}\cap\Omega_X$. Moreover, by Corollary VII.3.7(a) of \cite{kunen}, $\mathcal{R} = \mathcal{W}\bigcup \mathcal{D}$ is dense. By Theorem \ref{ellentuck}, $\mathcal{R}$, $\mathcal{W}$ and $\mathcal{D}$ are completely Ramsey. Thus, for the given $(s,B)\in\mathcal{M}_A$ there is a $B_1\in[B]^{\aleph_0}\cap\Omega_X$ such that $[s,B_1]\subset \mathcal{R}$, or $[s,B_1]\cap \mathcal{R}=\emptyset$; since $\mathcal{R}$ is dense and $[s,B_1]$ is nonempty and open we have $[s,B_1]\subset\mathcal{R}$. But now $\mathcal{W}$ is completely Ramsey and so there is a $B_2\in[B_1]^{\aleph_0}\cap\Omega_X$ with $[s,B_2]\subset \mathcal{W}$, or $[s,B_2]\cap\mathcal{W}=\emptyset$. Since $[s,B_2]\subset[s,B_1]$, we have that $[s,B_2]\subset \mathcal{W}$ or $[s,B_2]\subset\mathcal{D}$. In either case we have $(s,B_2)\forces \psi$, or $(s,B_2)\forces\neg\psi$.

\subsection{Proof of $(2)\Rightarrow (1)$:} This proof takes more work. We show that in fact $(2)$ implies that $\Omega_X\rightarrow(\Omega_X)^2_2$ holds. To see this, assume on the contrary that $\Omega_X\rightarrow(\Omega_X)^2_2$ fails. Choose a countable $A\in\Omega_X$ and a function $f:[A]^2\rightarrow \{0,1\}$ which witness this failure. Enumerate $A$ bijectively as $(a_n:n\in\naturals)$ and build the following corresponding partition tree:\\
$T_{\emptyset} = A$. $T_{(i)}:=\{a_n:n>1 \mbox{ and }f(\{a_1,a_n\}) = i\}$. For $\sigma\in^{<\omega}\{0,1\}$ of length $m$ for which $T_{\sigma}\in \Omega_X$, $T_{\sigma\frown(i)}:=\{a_n\in T_{\sigma}:n>m \mbox{ and }f(\{a_m,a_n\}) = i\}$. 

Observe that for each $\sigma$ with $T_{\sigma}\in\Omega_X$ we have $T_{\sigma\frown(0)}\in\Omega_X$ or $T_{\sigma\frown(1)}\in\Omega_X$. For each $n$, define $\mathcal{T}_n:=\{T_{\sigma}\in\Omega_X: length(\sigma)=n\}$. Then we have from the definitions that:
\begin{enumerate}
\item{For each $B\in[A]^{\aleph_0}\cap\Omega_X$ and for each $n$ there is a $T\in\mathcal{T}_n$ with $B\cap T\in\Omega_X$.}
\item{For each $n$, for each $T\in\mathcal{T}_{n+1}$ there is a unique $T^{\prime}\in\mathcal{T}_n$ with $T\subset T^{\prime}$.}
\item{For each $n$ and $\sigma$, if $a_n\in T_{\sigma}$, then $n> m = length(\sigma)$.} 
\end{enumerate}

{\flushleft{\bf Claim 1:}} If there is a $B\in[A]^{\aleph_0}\cap\Omega_X$ such that for each $n$ there is a $T\in\mathcal{T}_n$ with $B\setminus\{a_j:j\leq n\}\subset T$, then there is a $C\in[A]^{\aleph_0}\cap\Omega_X$ such that $f$ is constant on $[C]^2$.\\
For let such a $B$ be given. Since the elements of $\mathcal{T}_1$ are pairwise disjoint, choose the unique $i_1\in\{0,1\}$ with $B\setminus\{a_1\} \subset T_{(i_1)}$. Letting $T$ be the unique element of $\mathcal{T}_2$ with $B\setminus\{a_1,a_2\}\subset T$, we see that $T\subset T_{(i_1)}$, and so for a unique $i_2\in\{0,1\}$, $B\setminus\{a_1,a_2\}\subset T_{(i_1,i_2)}$. Arguing  like this we find an infinite sequence $(i_j:j<\infty)$ in $^{\naturals}\{0,1\}$ such that for each $m$, $B\setminus\{a_1,\cdots,a_m\}\subset T_{(i_1,\cdots,i_m)}$.

Write $B =\{a_{n_j}:j\in\naturals\}$ where $n_i<n_j$ whenever $i<j$. Put $B_1=\{a_{n_j}: i_{n_j}=1\}$ and $B_0=\{a_{n_j}: i_{n_j}=0\}$. Then $B_0\in\Omega_X$, or $B_1\in\Omega_X$. In the former case $f$ is constant of value $0$ on $[B_0]^2$, and in the latter case $f$ is constant of value $1$ on $[B_1]^2$. This completes the proof of Claim 1.

Note that the conclusion of Claim 1 holds also if instead we hypothesize that $B\in[A]^{\aleph_0}\cap\Omega_X$ is such that for each $n$ with $a_n\in B$ there is a $T\in \mathcal{T}_{n}$ with $B\setminus\{a_j:j\le n\}\subset T$.

Since we are assuming that there is no $B\in[A]^{\aleph_0}\cap\Omega_X$ with $f$ constant on $[B]^2$, we get: There is no $B\in[A]^{\aleph_0}\cap\Omega_X$  such that for each $n$ with $a_n\in B$ there is a $T\in\mathcal{T}_n$ with $B\setminus\{a_j:j\leq n\}\subset T$. Indeed, this is equivalent to: 
\begin{quote}
For each $B\in[A]^{\aleph_0}\cap\Omega_X$  there is an $n$ with $a_n\in B$ but for each $T\in\mathcal{T}_n$ we have $B\setminus\{a_j:j\leq n\}\not\subset T$.
\end{quote}

In what follows we will use ${\stackrel{\bullet}{a}}$ to denote the canonical name of the ground model object $a$ in the forcing language. 
Define the $\mathcal{M}_A$-name 
\[
  \Gamma:=\{({\stackrel{\bullet}{a}}_n,(s,B)): (s,B)\in \mathcal{M}_A \mbox{ and } a_n\in s\}.
\]
Then for each $\mathcal{M}_A$-generic filter $G$ we have 
\[
  \Gamma_G = \cup\{s\in[A]^{<\aleph_0}:(\exists B\in[A]^{\aleph_0}\cap\Omega_X)((s,B)\in G)\}.
\]

{\flushleft{\bf Claim 2:}} $(\emptyset,A)\forces ``(\exists n)(\forall T\in{\stackrel{\bullet}{\mathcal{T}}_n})(\Gamma\setminus\{\stackrel{\bullet}{a}_j:j\le n\}\not\subseteq T)"$.\\
For suppose that on the contrary $(s,B)\forces ``(\forall n)(\exists T\in{\stackrel{\bullet}{\mathcal{T}}}_n)(\Gamma\setminus\{\stackrel{\bullet}{a}_j:j\le n\}\subseteq T)"$. Since we have $B\in[A]^{\aleph_0}\cap\Omega_X$, choose an $n_1$ so that $B\setminus\{a_j:j\le n_1\}$ is not a subset of any $T\in\mathcal{T}_n$. Then choose a $T_{n_1}\in\mathcal{T}_{n_1}$ so that $B\cap T_{n_1}\in\Omega_X$. Also choose $a_m\in B\setminus(T_{n_1}\cup\{a_j:j\le n_1\})$. Put $B^{\prime} = B|\{a_j:j\le m_1\}$ and put $t = s\cup\{a_m\}$. Then as $(s,B)\forces ``\stackrel{\bullet}{\mathcal{T}}_{n_1}\mbox{ is a disjoint family}"$ and 
$(t, B^{\prime}\cap T_{n_1})\forces ``(\Gamma\setminus\{\stackrel{\bullet}{a}_j:j\le \stackrel{\bullet}{n}_1\})\cap \stackrel{\bullet}{T}_{n_1} \neq\emptyset"$.
\begin{equation}\label{force1}
   (t, B^{\prime}\cap T_{n_1})\forces ``\Gamma\setminus\{\stackrel{\bullet}{a}_j:j\le \stackrel{\bullet}{n}_1\}\subset\stackrel{\bullet}{T}_{n_1}".
\end{equation}
But evidently we also have 
\begin{equation}\label{force2}
   (t, B^{\prime}\cap T_{n_1})\forces ``\stackrel{\bullet}{a}_m\in(\Gamma\setminus\{\stackrel{\bullet}{a}_j:j\le \stackrel{\bullet}{n}_1\})\setminus \stackrel{\bullet}{T}_{n_1}".
\end{equation}
Thus we have a condition forcing contradictory statements, a contradiction. It follows that Claim 2 holds.

Now we construct a sentence $\Psi(\Gamma)$ in the forcing language:
\[
  ``\Gamma\cap\{\stackrel{\bullet}{a}_j:j<n\} \mbox{ is even for the least }n \mbox{ with } \stackrel{\bullet}{a}_n\in\Gamma \mbox{ and for all }T \in \stackrel{\bullet}{\mathcal{T}}_n \Gamma\setminus \{\stackrel{\bullet}{a}_j:j\le n\}\not\subset T,\, "
\]
By hypothesis 2 of the theorem, choose a $B\in[A]^{\aleph_0}\cap \Omega_X$ such that $(\emptyset,B)$ decides $\Psi(\Gamma)$. 

Choose $k_1$ minimal so that $a_{k_1}\in B$ and for each $T\in\mathcal{T}_{k_1}$ we have $B\setminus\{a_j:j\leq k_1\}\not\subseteq T$. 
Put $B_1 = B\setminus\{a_j:j\leq k_1\}$ and choose $T_{k_1}\in\mathcal{T}_{k_1}$ so that $C_1:=B_1\cap T_{k_1}\in\Omega_X$. Choose $\ell_1$ so that $a_{\ell_1}\in B_1\setminus T_{k_1}$. Then $(\{a_{k_1},a_{\ell_1}\}, C_1) < (\emptyset,B)$ and so also $(\{a_{k_1},a_{\ell_1}\}, C_1)$ decides $\Psi(\Gamma)$. 

By the construction of $C_1$ we see that for $T^{\prime}\in\mathcal{T}_{k_1}\setminus\{T_{k_1}\}$, 
also $(\{a_{k_1},a_{\ell_1}\}, C_1)\forces ``\Gamma\setminus\{\stackrel{\bullet}{a}_j:j\le \stackrel{\bullet}{k}_1\}\not\subseteq \stackrel{\bullet}{T^{\prime}}"$. 
And since $a_{\ell_1}\not\in T_{k_1}$ we also have 
$(\{a_{k_1},a_{\ell_1}\}, C_1)\forces ``\Gamma\setminus\{\stackrel{\bullet}{a}_j:j\le \stackrel{\bullet}{k}_1\}\not\subseteq \stackrel{\bullet}{T}_{k_1}"$. 
Moreover, $(\{a_{k_1},a_{\ell_1}\}, C_1)\forces ``\Gamma\cap\{\stackrel{\bullet}{a}_j:j< \stackrel{\bullet}{k}_1\}=\emptyset"$. 
Since $k_1$ was chosen minimal and $a_{k_1}\in B$, the least $n$ having the properties of $k_1$ is $k_1$.
It follows that $(\{a_{k_1},a_{\ell_1}\}, C_1)\forces \Psi(\Gamma)$, and as $(\emptyset,B)$ already decides $\Psi(\Gamma)$, we have 
\begin{equation}\label{positiveforce}
(\emptyset,B)\forces \Psi(\Gamma).
\end{equation}

Now repeat the previous construction starting with $C_1$ in place of $B$. Choose $k_2$ minimal so that $a_{k_2}\in C_1$ and for each $T\in\mathcal{T}_{k_2}$ we have $C_1\setminus\{a_j:j\leq k_2\}\not\subseteq T$. Since $a_{k_1}\not\in C_1$, we have $k_2>k_1$. Put $B_2 = C_1\setminus\{a_j:j\leq k_2\}$ and choose $T_{k_2}\in\mathcal{T}_{k_2}$ so that $C_2:=B_2\cap T_{k_2}\in\Omega_X$. Choose $\ell_2$ so that $a_{\ell_2}\in B_2\setminus T_{k_2}$. Then $(\{a_{k_1}, a_{k_2}, a_{\ell_2}\}, C_2) < (\emptyset,B)$ and so also $(\{a_{k_1}, a_{k_2}, a_{\ell_2}\}, C_2)$ decides $\Psi(\Gamma)$. 
By the construction of $C_2$ we see that for $T^{\prime}\in\mathcal{T}_{k_2}\setminus\{T_{k_2}\}$, 
also $(\{a_{k_1}, a_{k_2}, a_{\ell_2}\}, C_2)\forces ``\Gamma\setminus\{\stackrel{\bullet}{a}_j:j\le \stackrel{\bullet}{k}_2\}\not\subseteq \stackrel{\bullet}{T^{\prime}}"$. 
And since $a_{\ell_2}\not\in T_{k_2}$ we also have 
$(\{a_{k_1}, a_{k_2}, a_{\ell_2}\}, C_2)\forces ``\Gamma\setminus\{\stackrel{\bullet}{a}_j:j\le \stackrel{\bullet}{k}_2\}\not\subseteq \stackrel{\bullet}{T}_{k_2}"$. 
By minimality of $k_2$ and the fact that $a_{k_2}\in C_2$, we get that the minimal $n$ with these properties of $k_2$ is $k_2$:
However, $(\{a_{k_1}, a_{k_2}, a_{\ell_2}\}, C_2)\forces ``\Gamma\cap \{\stackrel{\bullet}{a}_j:j\le \stackrel{\bullet}{k}_2\}=\{\stackrel{\bullet}{a}_{k_1}\}"$. This means that $(\{a_{k_1}, a_{k_2}, a_{\ell_2}\}, C_2)\forces \neg\Psi(\Gamma)$. Since $(\{a_{k_1}, a_{k_2}, a_{\ell_2}\}, C_2) < (\emptyset,B)$ and $(\emptyset,B)$ already decides $\Psi(\Gamma)$, we find that 

\begin{equation}\label{negativeforce}
(\emptyset,B)\forces\neg\Psi(\Gamma).
\end{equation}
Since (\ref{positiveforce}) and (\ref{negativeforce}) yield a contradiction, the hypothesis that $\Omega\rightarrow(\Omega)^2_2$ fails is false. This completes the proof of $(2)\Rightarrow(1)$ of Theorem \ref{truthlemma}. $\epf$ 

{\flushleft{\bf Remarks:}} The above result is again given for $\Omega$, but a study of the proofs will reveal that these equivalences hold for several other families $\mathcal{A}$, including the examples mentioned earlier. Theorem \ref{truthlemma} has several consequences that will be explored elsewhere. One of the mentionable consequences is that forcing with $\mathcal{M}_A$ preserves cardinals, and in the generic extension the only groundmodel sets of reals having $\sone(\Omega,\Omega)$ are the countable sets. And a countable support iteration of length $\aleph_2$ over a ground model satisfying the Continuum Hypothesis gives a model of Borel's Conjecture, just like the usual Mathias reals iteration does - \cite{JEB}. 

In closing: Analogous results can be proved for the selection principle $\sfin(\mathcal{A},\mathcal{A})$ and its relatives. These will be reported elsewhere.

\begin{center}{\bf Appendix A: Regarding Theorem \ref{soneramsey}:}
\end{center}

Strictly speaking, the only equivalence that has been explicitly proved in the literature is the equivalence of $\sone(\Omega,\Omega)$ with $\Omega\longrightarrow(\Omega)^2_2$, with the remark that the techniques used to prove this case yield by an induction the full version that for all finite $n$ and $k$ we have $\Omega\longrightarrow(\Omega)^n_k$. It is perhaps worth putting down the main elements of such an argument explicitly for future reference. The only implication we need to prove is the implication that $\Omega\longrightarrow(\Omega)^2_2$ implies that for all $n$ and $k$ $\Omega\longrightarrow(\Omega)^n_k$.
{\flushleft{\bf Claim 1:}} $\Omega\longrightarrow(\Omega)^2_2$ implies that for each $k>1$, $\Omega\longrightarrow(\Omega)^2_k$.\\ 
This can be done by induction on $k+1$. For $k=1$ this is the hypothesis. Assuming we have proven the implication for $j\le k$, consider a countable $\omega$-cover $\mathcal{U}$ of $X$ and a coloring $f:[\mathcal{U}]^2\longrightarrow\{1,\cdots,k+1,k+2\}$.  Define a new coloring $g$ so that
\[
  g(\{U,V\}) = \left\{\begin{tabular}{ll}
                      f(\{U,V\}) & \mbox{ if }f(\{U,V\})$<$ k+1.\\
                      k+1        & \mbox{ otherwise}
                      \end{tabular}\right. 
\]
Applying the induction hypothesis we find an $\omega$-cover $\mathcal{V}\subset\mathcal{U}$ and an $i\in\{1,\cdots,k+1\}$ such that $g(\{U,V\})=i$ for all $\{U,V\}\in\mathcal{V}^2$. If $i<k+1$ then indeed $\mathcal{V}$ works for $f$. Else, $\mathcal{V}$ is an $\omega$-cover on whose pairs $f$ takes values $k+1$ or $k+2$, and now apply $\Omega\longrightarrow(\Omega)^2_2$.\\
{\flushleft{\bf Claim 2:}} For $n>2$ and $k>1$, $\Omega\longrightarrow(\Omega)^n_k$ implies $\Omega\longrightarrow(\Omega)^2_k$.\\
This can be done by starting with a countable $\omega$-cover $\mathcal{U}$ and a coloring $f:[\mathcal{U}]^2\longrightarrow\{1,\cdots,k\}$. Enumerate $\mathcal{U}$ bijectively as $\{U_m:m\in\naturals\}$. Define $g:[\mathcal{U}\}]^n\longrightarrow\{1,\cdots,k\}$ by 
\[ 
  g(\{U_{i_1}\cdots,U_{i_n}\}) = f(\{U_{i_1},U_{i_2}\}),
\]
where we list the $n$-tuples according to increasing index in the chosen enumeration. Apply $\Omega\longrightarrow(\Omega)^n_k$.\\
{\flushleft{\bf Claim 3:}} For $n>1$ and $k>1$, $\Omega\longrightarrow(\Omega)^n_k$ implies $\Omega\longrightarrow(\Omega)^{n+1}_k$.\\
To prove this we use the fact that For $n>1$ and $k>1$, $\Omega\longrightarrow(\Omega)^n_k$ implies $\Omega\longrightarrow(\Omega)^2_2$, which in turn implies that ONE has no winning strategy in the game $\gone(\Omega,\Omega)$. 

Let a countable $\omega$-cover $\mathcal{U}$ be given, as well as $f:[\mathcal{U}]^{n+1}\longrightarrow\{1,\cdots,k\}$. Enumerate $\mathcal{U}$ bijectively as $\{U_m:m\in\naturals\}$. Define a strategy $F$ for ONE in the game $\gone(\Omega,\Omega)$ as follows:\\
Fix $U_1$ and define 
\[
  g_1:[\mathcal{U}\setminus\{U_1\}]^n\longrightarrow\{1,\cdots,k\}
\]
by $g_1(\mathcal{V}) = f(\{U_1\}\bigcup\mathcal{V})$. Using $\Omega\longrightarrow(\Omega)^n_k$, fix an $i_1\in\{1,\cdots,k\}$ and an $\omega$-cover $\mathcal{U}_1\subset \mathcal{U}$ such that $g_1(\mathcal{V})=i_1$ for each $\mathcal{V}\in[\mathcal{U}_1]^n$. Declare ONE's move to be $F(\emptyset) = \mathcal{U}_1$. 

When TWO responds with $T_1 = U_{n_1}\in F(\emptyset)$, ONE first defines
\[
  g_2:[\mathcal{U}_1\setminus\{U_j:j\le n_1\}]^n\longrightarrow\{1,\cdots,k\}
\]
by $g_2(\mathcal{V}) = f(\{U_{n_1}\}\bigcup\mathcal{V})$. Then, using $\Omega\longrightarrow(\Omega)^n_k$, fix an $i_{n_1}\in\{1,\cdots,k\}$ and an $\omega$-cover $\mathcal{U}_2\subset \mathcal{U}_1\setminus\{U_j:j\le n_1\}$ such that $g_2(\mathcal{V})=i_{n_1}$ for each $\mathcal{V}\in[\mathcal{U}_2]^n$. Declare ONE's move to be $F(T_1) = \mathcal{U}_2$.

When TWO responds with $T_2 = U_{n_2}\in F(T_1)$, ONE first defines
\[
  g_3:[\mathcal{U}_2\setminus\{U_j:j\le n_2\}]^n\longrightarrow\{1,\cdots,k\}
\]
by $g_3(\mathcal{V}) = f(\{U_{n_2}\}\bigcup\mathcal{V})$. Then, using $\Omega\longrightarrow(\Omega)^n_k$, fix an $i_{n_2}\in\{1,\cdots,k\}$ and an $\omega$-cover $\mathcal{U}_3\subset \mathcal{U}_2\setminus\{U_j:j\le n_2\}$ such that $g_3(\mathcal{V})=i_{n_2}$ for each $\mathcal{V}\in[\mathcal{U}_3]^n$. Declare ONE's move to be $F(T_1,T_2) = \mathcal{U}_3$.

This describes ONE's strategy in this game. Since it is not winning for ONE, we find a play $F(\emptyset), T_1, F(T_1), T_2, F(T_1,T_2), T_3, \dots$ which is lost by ONE. Associated with this play we have an increasing infinite sequence $n_1< n_2<\cdots <n_k<\cdots$ for which $T_k = U_{n_k}$, all $k$, and a sequence $i_{n_k}, k\in \naturals$ of elements of $\{1,\cdots,k\}$, and a sequence $\mathcal{U}_n,\, n\in\naturals$, of $\omega$-covers such that:
\begin{enumerate}
\item{For each $m$, $T_m=U_{n_m}\in\mathcal{U}_m\subset\mathcal{U}_{m-1}\setminus\{U_{n_j}:j\le m-1\}$.}
\item{For each $m$, $f(\{T_m\}\bigcup\mathcal{V}) = i_{n_m}$ whenever $\mathcal{V}\in[\mathcal{U}_{m+1}]^n$.}
\item{$\{T_m:m\in\naturals\}\subset\mathcal{U}$ is an $\omega$-cover.} 
\end{enumerate}
Fix an $i$ such that $\mathcal{W} = \{T_m:i_{n_m}=i \mbox{ and }m>n\}$ is an $\omega$-cover. Then for each $\mathcal{V}\in[\mathcal{W}]^{n+1}$ we have $f(\mathcal{V})=i$. 

\begin{center}{\bf Acknowledgements}\end{center}
I thank the referee for a careful reading of the paper, and for very useful suggestions. I also thank the Boise Set Theory seminar for useful suggestions and remarks during a series of lectures on these results. And finally, I thank the organizing committee of the third Workshop on Coverings, Selections and Games in Topology, hosted in April 2007 in Vrnjacka Banja, Serbia, for the opportunity to present some of these results at the workshop.

Address:\\
Department of Mathematics\\
Boise State University\\
Boise, ID 83725.\\
e-mail: marion@math.boisestate.edu

\begin{thebibliography}{}
\bibitem{JEB} J.E. Baumgartner, \emph{Iterated forcing}, {\bf London Mathematical Society Lecture Note Series} 87 (1983), 1 - 59.
\bibitem{dkm} G. di Maio, Lj.D.R. Ko\v{c}inac and E. Meccariello, \emph{Applications of k-covers}, {\bf Acta Mathematica Sinica} (English Series) 22 (2006),  1151 - 1160.
\bibitem{El} E. Ellentuck, \emph{A new proof that analytic sets are Ramsey}, {\bf The Journal of Symbolic Logic} 39:1 (1974), 163 - 165.
\bibitem{Galvinnotices} F. Galvin, \emph{A generalization of Ramsey's Theorem}, {\bf Notices of the American Mathematical Society} 15 (1968), p. 548. Abstract 68T-368.
\bibitem{G-P} F. Galvin and K. Prikry, \emph{Borel sets and Ramsey's Theorem}, {\bf The Journal of Symbolic Logic} 38:2 (1973), 193 - 198.
\bibitem{GN} J. Gerlits and Zs. Nagy, \emph{ Some properties of C(X), I}, {\bf Topology and its Applications}    14 (1982), 151 -- 161.
\bibitem{coc2} W. Just, A.W. Miller, M. Scheepers and P.J. Szeptycki, {\em Combinatorics of open covers (II)}, {\bf Topology and its Applications} 73 (1996), 241 - 266.
\bibitem{kunen} K. Kunen, \emph{Set Theory: An introduction to independence proofs}, {\bf Studies in Logic and the Foundations of Mathematics} 102 (1980).
\bibitem{Ro38} F. Rothberger, {\em Eine Versch\"arfung der Eigenschaft {\sf C}}, {\bf Fundamenta Mathematicae} 30 (1938), 50 - 55.
\bibitem{Sakai} M. Sakai, {\em Property {\sf C}'' and function spaces}, {\bf Proceedings of the American Mathematical Society} 104 (1988), 917 - 919.
\bibitem{ST} N. Samet and B. Tsaban, \emph{Ramsey theory of open covers: Lecture 3}, preprint. 
\bibitem{MSBct}  M.  Scheepers, \emph{The least cardinal for which the Baire category theorem fails}, {\bf Proceedings of the American Mathematical Society} 125:2 (1997), 579 - 585. 
\bibitem{coc1} M. Scheepers, {\em Combinatorics of open covers (I): Ramsey Theory}, {\bf Topology and its Applications} 69 (1996), 31-62.
\bibitem{coc3} M. Scheepers, {\em Combinatorics of open covers (III): games, ${\sf C}_p(X)$}, {\bf Fundamenta Mathematicae} 152 (1997), 231 - 254.
\bibitem{coc5} M. Scheepers, \emph{Combinatorics of open covers (V): Pixley-Roy spaces of sets of reals and $\omega$-covers}, {\bf Topology and its Applications}102 (2000), 13 - 31.
\bibitem{MSDL} M. Scheepers, \emph{$\sone(\mathcal{A},\mathcal{B})$ in distributive lattices}, {\bf Quaderni Mathematicae} (to appear)
\end{thebibliography}
\end{document}